\documentclass{amsart}
\usepackage{amsmath,amsfonts,amscd,latexsym,amssymb,epsfig,float,
afterpage}

\newtheorem{theorem}{Theorem}[section]
\newtheorem{lemma}[theorem]{Lemma}
\newtheorem{proposition}[theorem]{Proposition}

\theoremstyle{definition}
\newtheorem{definition}[theorem]{Definition}

\newtheorem{corollary}[theorem]{Corollary}
\newtheorem{remark}[theorem]{Remark}

\newtheorem{proposition+definition}[theorem]{Proposition \& Definition}

\theoremstyle{remark}

\numberwithin{equation}{section}

\newcommand{\cpxm}{{\mathbb{C}}}
\newcommand{\ratm}{{\mathbb{Q}}}
\newcommand{\itgm}{{\mathbb{Z}}}

\newcommand{\Hom}{\operatorname{Hom}}

\newcommand{\Res}{\operatorname{Res}}

\newcommand{\Extmhs}{\operatorname{Ext_{{\scriptscriptstyle MHS}}}}

\newcommand{\Pic}{\operatorname{Pic}}
\newcommand{\Jac}{\operatorname{Jac}}
\newcommand{\Period}{\operatorname{\Pi}}
\newcommand{\Mat}{\operatorname{Mat}}
\newcommand{\tr}{\operatorname{tr}}
\newcommand{\diag}{\operatorname{diag}}

\newcommand{\Gr}{\operatorname{Gr}}

\newcommand{\topa}[2]{\genfrac{}{}{0pt}{}{#1}{#2}}
\begin{document}
\renewcommand{\labelenumi}{(\roman{enumi})}

\title[The MHS on the Fundamental Group of a Punctured
Riemann Surface]{The Mixed Hodge Structure on the 
Fundamental Group of a Punctured Riemann Surface}

\author{Rainer H.~Kaenders}

\address{Mathematisches Institut, Heinrich-Heine-Universit{\"a}t 
D{\"u}sseldorf, Universit{\"a}tsstra{\ss}e, 40225 D{\"usseldorf}, Germany} 
\email{kaenders@cs.uni-duesseldorf.de}

\thanks{partly supported by grant ERBCHBICT930403 (HCM) from the
Eurpean Community and The Netherlands Organisation for Scientific
Research (NWO)}

\subjclass{Primary 14H40, 14H30; Secondary 14F35 }

\date{\today}

\begin{abstract}
Given a compact Riemann surface $\bar{X}$ of genus $g$ and a point $q$ on 
$\bar{X}$, we consider $X:=\bar{X}\setminus\{q\}$ with a basepoint $p\in X$. 
The extension of mixed Hodge structures, given by the weights -1 and -2, 
of the mixed Hodge 
structure on the fundamental group (in the sense of Hain) is studied.
We show that it naturally corresponds on the one hand to the element 
$(2g\,q-2\,p-K)$ in $\Pic^0(\bar{X})$, where $K$ represents the canonical 
divisor, and on the other hand to the respective extension of 
$\bar{X}$. Finally, we prove a pointed Torelli theorem for punctured 
Riemann surfaces.
\end{abstract}

\maketitle

% Der folgende Befehl laesst nur die "Uberschriften bis zur zweiten Stufe
% im Inhaltsverzeichnis erscheinen.
 
\setcounter{tocdepth}{1}

% \tableofcontents

\section*{Introduction}

Let $q$ be a point in a compact Riemann surface $\bar{X}$ of genus $g$. 
In this paper we
want to study the complement of $q$ in $\bar{X}$, i.~e.~$X:=
\bar{X}\setminus \{q\}$ with a basepoint $p$.
We refer to this situation as to a {\it punctured Riemann surface 
$X$ with puncture $q\in\bar{X}$ and basepoint $p$}.

For compact Riemann surfaces, Hain and Pulte
(\cite{Hain-the-geometry}, \cite{Pulte}) proved that the extension of
mixed Hodge structure, which is given by the weights -1 and -2 of the 
mixed Hodge structure (MHS) on the fundamental
group, determines the base point (see Theorem \ref{injection}).
From this result and from the classical Torelli 
theorem they derived a {\it
pointed Torelli theorem} (see Theorem \ref{pointed Torelli theorem}) 
as a corollary.

For a punctured Riemann surface $(\bar{X}\setminus\{q\},p)$, 
the extension $w_{pq}$
given by the weights -1 and -2 of the MHS on the fundamental group is one 
dimension bigger than in the compact case. 
We show that it naturally corresponds on the one hand to the element
\vspace{-3ex}
\[
(2g\,q-2\,p-K) \text{  in  } \Pic^0(\bar{X}),
\]
where $K$ represents the canonical
divisor, and on the other hand to the respective extension of
$\bar{X}$ (see Theorem \ref{Hauptsatz}).
This may have implications on possible normal functions on the moduli space
of complex projective curves (cf.~\cite{Hain-Looijenga}, 7.4).

Finally we prove that this extension $w_{pq}$ determines both, the
basepoint $p$ and the `removed point' $q$. This, together with the pointed
Torelli theorem of Hain and Pulte yields a
{\it punctured pointed Torelli theorem} 
(see Theorem \ref{two-pointed Torelli theorem}).

{\it Acknowledgements:}
The contents of this paper were part of my PhD-thesis which was accepted 
at the Catholic University of Nijmegen, NL in october 1997. 
I am indepted to my advisor Joseph Steenbrink, who posed the question 
for the above extension of MHS, and to Richard Hain and Eduard 
Looijenga for interesting discussions on the subject.

\section{Extensions and the Theta Divisor} \label{Haensel}

For the definition of iterated integrals and of the MHS on the 
fundamental group we refer to the introductory article \cite{Hain-the-geometry}.

Let $\bar{X}$ be a compact Riemann surface of genus $g$ and let $q$ be a point
on $\bar{X}$. We consider the pointed space $(X,p)$, where 
$X:=\bar{X}\setminus\{q\}$ and $p$ is a basepoint on $X$.
Denote by $J\subset \itgm\pi_1(X,p)$ and $\bar{J}\subset \itgm\pi_1(\bar{X},p)$
the augmentation ideals in the group rings of the respective fundamental 
groups. Note that ${J}\big/{J}^{2}=H_1(X)$ 
resp.~$\bar{J}\big/\bar{J}^{2}=H_1(\bar{X})$, and since we remove only a single 
point from $\bar{X}$, we have that 
$X\hookrightarrow\bar{X}$ induces an isomorphism of pure Hodge structures
between $H_1(X)$ and $H_1(\bar{X})$, both of weight -1. This allows us to 
identify these two Hodge structures and their duals. We will write just
$H_1$ and $H^1$. 
The {\it MHS on the fundamental group} $\pi_1(X,p)$ 
resp.~$\pi_1(\bar{X},p)$ consists by definition of MHSs on the integral 
lattices ${J}\big/{J}^{s+1}$ resp.~$\bar{J}\big/\bar{J}^{s+1}$ for $s\ge 2$.

This definition of the MHSs is possible because of Chen's 
$\pi_1$-De Rham-Theorem, telling us that integration of iterated integrals 
yields isomorphisms
\begin{equation} \label{Goldesel}
\begin{split} 
H^0\bar{B}_s\big(E^\bullet(\bar{X}\log q),\; p \big)
& \xrightarrow{\cong} \Hom_\itgm(J\big/J^{s+1},\;\cpxm)=: 
\left( J\big/J^{s+1} \right)_\cpxm^*\qquad \text{ resp. } \\
H^0\bar{B}_s\big(E^\bullet(\bar{X}),\; p \big)
& \xrightarrow{\cong} \Hom_\itgm(\bar{J}\big/\bar{J}^{s+1},\;\cpxm)=:
\left( \bar{J}\big/\bar{J}^{s+1} \right)_\cpxm^*.
\end{split}
\end{equation}
Here $E^\bullet(\bar{X} \log q)$ denotes the differential graded algebra (dga)
of $C^\infty$-forms on $X=\bar{X}\setminus\{q\}$ with logarithmic 
singularities at $q$ and 
$E^\bullet(\bar{X})$ denotes the dga of smooth 
complex valued forms on $\bar{X}$.
The objects on the left of \eqref{Goldesel} are the complex 
vector spaces of iterated integrals
of length $\le s$, which are homotopy functionals -- considered as
functions on the fundamental group. These vector spaces can be described
purely algebraically in terms of the augmented dga's 
$E^\bullet(\bar{X}\log q)$ and $E^\bullet(\bar{X})$. This 
is part of a general construction, {\it the reduced bar construction},
whence the elaborate notation (cf.~ \cite{Chen-RBC} or
\cite{Hain-de-Rham-homotopy}). Here we identify these different descriptions 
deliberately.

In the two cases under consideration, the weight filtration $W_{\bullet}$
is already given on the lattices 
${J}\big/{J}^{s+1}$ resp.~$\bar{J}\big/\bar{J}^{s+1}$ by 
the $J$-adic filtration, i.~e.~
\[
W_{-l} \; {J}\big/{J^{s+1}}={J^{l}}\big/{J^{s+1}} 
\text{  resp.  } W_{-l} \; {\bar{J}}\big/{\bar{J}^{s+1}}=
{\bar{J}^{l}}\big/{\bar{J}^{s+1}} 
\text{  for  } l>0.
\] 
For $l=1$ and $s=1$ we recover the pure Hodge structure on homology, 
i.~e.~ $W_{-1} {J}\big/{J^{2}}$ = ${J}\big/{J^{2}}$ = $H_1$ =
${\bar{J}}\big/{\bar{J}^{2}}$ = $W_{-1}\bar{J}\big/{\bar{J}^{2}}$.
The weights -1 and -2 give rise to two extensions of MHSs $w_{pq}$ and $w_p$,
related by the following commutative diagram
\begin{equation}\label{Hexe}
\begin{CD}
0 @>>> {J^2}\big/{J^{3}} @>>> {J}\big/{J^{3}} @>>> {J}\big/{J^{2}} @>>> 0 & 
\qquad ;\; w_{pq} \\
& & @VVV @VVV @VV=V & \\
0 @>>> \bar{J}^2\big/\bar{J}^{3} @>>> \bar{J}\big/\bar{J}^{3} @>>> 
\bar{J}\big/\bar{J}^{2} @>>> 0  &\qquad ;\; w_p. 
\end{CD}
\end{equation}

The multiplication in the group rings defines surjective maps
$J\big/ J^2\otimes J\big/ J^2\rightarrow J^2\big/ J^3$ and 
$\bar{J}\big/ \bar{J}^2\otimes \bar{J}\big/ \bar{J}^2\rightarrow 
\bar{J}^2\big/ \bar{J}^3$ whose dual morphisms are inclusions 
$(J^2\big/ J^3)^*\hookrightarrow H^1\otimes H^1$ and $(\bar{J}^2\big/ 
\bar{J}^3)^*\hookrightarrow H^1\otimes H^1$.
It is well-known (cf.~\cite{Hain-the-geometry}) that in both cases, the image 
of the above inclusions coincides with the kernel of the cup-product.
Hence the inclusions give isomorphisms $(J^2\big/ J^3)^*\cong H^1\otimes H^1$
and $(\bar{J}^2\big/ \bar{J}^3)^*\cong K$, where $K:=\ker\{
\cup:H^1(\bar{X})\otimes H^1(\bar{X})\rightarrow H^2(\bar{X})\}$. As
$\cup$ is a morphism of Hodge structures, $K$ 
inherits a pure Hodge structure of weight 2 from $H^1\otimes H^1$. 
Dualizing the diagram \eqref{Hexe}
yields extensions of MHSs $m_{pq}$ and $m_p$ -- dual to $w_{pq}$ and $w_p$
\begin{equation*}
\begin{CD}
0 @>>> H^1 @>>> \left({J}\big/{J^{3}}\right)^* @>>> H^1\otimes H^1 @>>> 0 &
\qquad ;\; m_{pq} \\
& & @A=AA @AAA @AAA & \\
0 @>>> H^1 @>>> \left(\bar{J}\big/\bar{J}^{3}\right)^* @>>>
K @>>> 0  &\qquad ;\; m_p.
\end{CD}
\end{equation*}
Since $K$ is not a direct summand of $H^1\otimes H^1$ over $\itgm$,  
let us first clarify the nature of the embedding $K\hookrightarrow 
H^1\otimes H^1$.
Identify $H^2(\bar{X},\itgm)$ with $\itgm$. There is a bilinear form
\vspace{-2ex}
\[
b:\big(H^1_\itgm\otimes H^1_\itgm\big) 
\times\big(H^1_\itgm\otimes H^1_\itgm\big)
\longrightarrow\itgm,
\]
given by $b\left((x_1\otimes x_2),\,(y_1\otimes y_2)\right)
:= (x_1\cup y_2)\cdot(y_1\cup x_2)$, which has mixed signature and is
non degenerate. 
Consider the rank 1 sublattice of $H^1_\itgm\otimes H^1_\itgm$ orthogonal
to the kernel of the cup-product $K_\itgm\subset H^1_\itgm\otimes H^1_\itgm$ 
with respect to $b$.  
The following commutative diagram describes how $Q_\itgm$ is related to 
$K_\itgm\subset H^1_\itgm\otimes H^1_\itgm$ and $H^2(\bar{X},\itgm)$.
\vspace{-2ex}
\begin{equation*}
\begin{CD}
Q_\itgm\oplus K_\itgm  & \hookrightarrow & H^1_\itgm\otimes H^1_\itgm  \\
@V{\frac{1}{2g}\cup}VV  @VV{\cup}V \\
H^2(\bar{X},\itgm) @>{\cdot 2g}>> H^2(\bar{X},\itgm).
\end{CD}
\end{equation*}
$Q_\itgm$ is generated by one element ${\mathfrak X}$ in 
$H^1_\itgm\otimes H^1_\itgm$, which is
invariant under complex conjugation. Hence $H^1_\itgm\otimes H^1_\itgm$ 
induces on 
$Q_\itgm$ a $\itgm$-HS, isomorphic to $H^2(\bar{X},\itgm)$ or $\itgm(-1)$.
Note that over the rationals holds:
$K_\ratm\oplus Q_\ratm= H^1_\ratm\otimes H^1_\ratm$.

\begin{definition} \rm
Define $k_{pq}=
[0\rightarrow H^1 \rightarrow E_{pq} \rightarrow Q \rightarrow 0]
\in \Extmhs\left(Q;H^1\right)$ to be the restriction of the
extension $m_{pq}$ to an extension of $Q$ by $H^1$. 
\end{definition}

Let $\Psi:\Extmhs\left(Q;H^1\right)\xrightarrow{\cong}
\Pic^0(\bar{X})$ be the natural isomorphism (see \cite{Carlson}), then 
the main theorem of this paper is
\begin{theorem} \label{Hauptsatz}
\qquad $\Psi(k_{pq})=(2g\, q -2\, p -K) \quad\in \Pic^0(\bar{X}).$
\end{theorem}

\subsection{Riemann's Constant}

Let $u:\Pic^0 (\bar{X})\rightarrow \Jac(\bar{X})$ be the Abel-Jacobi map 
and define the divisor 
\[ {\textstyle
W_{p,g-1}:=\left\{\sum\nolimits_{j=1}^{g-1} u(q_j-p)\;\big| \; 
\sum\nolimits_{j=1}^{g-1} q_j\in \bar{X}^{(g-1)}\right\}.} 
\]
Denote the {\it theta divisor} on $\Jac(\bar{X})$ by $\Theta$ and
{\it Riemann's constant} by $\kappa_p\in\Jac(\bar{X})$, 
such that Riemann's classical 
theorem\footnote{Proofs of this theorem can be found in \cite{Riemann}
(VI, 22., pp.~132-136; XI, pp.~213-224) or \cite{Landfriedt}.
For proofs in modern language we refer to \cite{Lewittes}, 
\cite{Tata-on-theta-I}
(Theorem 3.1, pp.~149-151) or to \cite{Griffiths-Harris}.
In the theory of $\theta$-functions it is more convenient to
define $\varkappa_p$ in $\cpxm^g$ like \cite{Riemann}, \cite{Landfriedt},
\cite{Lewittes}, \cite{Fay}
(here Riemann's constant is defined as $-\varkappa_p$)
and \cite{Tata-on-theta-I} do.} 
reads: $\Theta=W_{p,g-1}+\kappa_p$.

Using the Riemann-Roch theorem one can prove that Riemann's constant 
$\kappa_p$
is related to the canonical divisor by the fact that for any divisor
$K$ of a holomorphic 1-form holds $u\big( (2g-2)p-K\big)=2\kappa_p$ and
that the canonical divisor is characterized by this equation (for a proof 
we refer to \cite{Griffiths-Harris}, p.~340).
Theorem \ref{Hauptsatz} is then a consequence of 
the following theorem, whose proof will be given in the sequel.

\begin{theorem}\label{Hans guck in die Luft}
\qquad $u\big(\Psi(k_{pq})+ 2g(p-q)\big)=2\kappa_p$.
\end{theorem}

The rest of this section is devoted to the proof of 
Theorem \ref{Hans guck in die Luft}. First 
we interpret the right hand side of the equation by means of an 
expression for $\kappa_p$
in terms of iterated integrals, as it was already known to Riemann.

To present this formula we need some more notation.
Denote by ${\gamma}_1,\ldots,{\gamma}_{2g}$
and ${\delta}$ a system of elements in $\pi_1(X,p)$
having the property, that
the fundamental group $\pi_1({X},p)$ is the quotient of the
free group $F\langle \gamma_1,\ldots,\gamma_{2g},\,\delta\rangle$
generated by the $\gamma_i$ and $\delta$ subject to the commutator relation 
\begin{equation} \label{Hannelore}
[\gamma_1,\; \gamma_{g+1}] \cdots
[\gamma_g,\; \gamma_{2g}]=\delta.
\end{equation}

Let $dz_1,\ldots,dz_g$ be a basis of holomorphic 1-forms on $\bar{X}$, such
that $\int_{\gamma_\nu} dz_i=\delta_{i\nu}$, i.~e.~the period matrix 
can be written $\Omega=(\omega_{i\mu})_{\topa{i=1,\ldots,g}{\mu=
1,\ldots,2g}}=\left(\Omega_1,\Omega_2\right) =\left(I,Z\right)$.
By Riemann's bilinear relations, $Z$ is a symmetric $g\times g$-matrix
with positive definite imaginary part. Having made these choices 
we may represent the Jacobian of $\bar{X}$ as
$\Jac (X):=\cpxm^g\big/ \Omega\itgm^{2g}$.
The following expresssion of $\kappa_p$ in terms of iterated integrals  
\begin{equation} \label{Rumpelstielzchen}
\kappa_p=\left[
-\sum_{\nu=1}^g \int_{\gamma_\nu} dz_i dz_\nu
+\frac{1}{2} \int_{\gamma_{g+i}} dz_i \right]_{i=1,\ldots,g} 
\in \Jac (\bar{X})
\end{equation}
was known to Bernhard Riemann in 1865 (see \cite{Riemann}, p.~213
or \cite{Fay}).

\subsection{Extension Data}

According to \cite{Carlson} we need two things for the computation of 
the extension data $k_{pq}\in \Extmhs(Q,H^1)$: a Hodge filtration 
preserving section 
$s_F:(Q,F^\bullet)\rightarrow (E_{pq},F^\bullet)$ and an integral retraction
$r_\itgm: (E_{pq})_\itgm \rightarrow H^1_\itgm$. 

Let $dx_1,\ldots,dx_{2g}$
be the real harmonic 1-forms such that $\int_{\gamma_j} dx_i = \delta_{ij}$.
Then a generator ${\mathfrak X}$ of $Q_\itgm$
is given by ${\mathfrak X}:=
\sum_{\nu=1}^g \left( [dx_\nu]\otimes[dx_{g+\nu}]
- [dx_{g+\nu}]\otimes[dx_\nu]\right) $.
Riemann's first bilinear relation tells us that ${\mathfrak X}\in 
F^1(H^1_\cpxm\otimes H^1_\cpxm)$. To be more precise, we can write 
\vspace{-2ex}
\begin{equation*} \label{Tina}
\sum\nolimits_{\nu=1}^g \left(dx_\nu\otimes dx_{g+\nu}
- dx_{g+\nu}\otimes dx_\nu \right)= \sum\nolimits_{j,k=1}^g
\left(a_{jk} dz_j\otimes d\bar{z}_k+ \bar{a}_{jk} d\bar{z}_j\otimes dz_k\right)
\end{equation*}
with $A=(a_{jk})_{jk}=\left(\bar{\Omega}_2\Omega_1^t-
\bar{\Omega}_1\Omega_2^t\right)^{-1}=(\bar{Z}-Z)^{-1}$. Observe:
$A^t=-\bar{A}$.

Define $\wedge{\mathfrak X}:=\sum_{\nu=1}^g \left(dx_\nu\wedge dx_{g+\nu}
- dx_{g+\nu}\wedge dx_\nu \right)\in F^1 E^2(\bar{X})$. The strictness of the
differential with respect to the Hodge filtration on 
$E^\bullet(\bar{X} \log q)$ implies that there is a $\mu_q\in 
F^1 E^1(\bar{X} \log q)$ such that $\wedge{\mathfrak X}+d\mu_q=0$.
This condition implies that
the iterated integral $\int {\mathfrak X}+\mu_q:=\int\sum_{\nu=1}^g
\left(dx_\nu dx_{g+\nu} - dx_{g+\nu} dx_\nu \right)+\mu_q$                      is a homotopy functional.    

A Hodge fitration preserving section $s_F$ is then defined by 
$s_F({\mathfrak X})=\int {\mathfrak X}+\mu_q$ and an integral retraction 
$r_\itgm$ is given by the map, which sends an iterated integral $\int I$
of length $\le 2$ with values in $\itgm$ to $r_\itgm(\int I):=
\sum_{j=1}^{2g}(\int_{\gamma_j} I)[dx_j]$.
Again a standard computation shows  
\[
\boxed{u\circ\Psi(k_{pq})=\left(\sum_{\nu=1}^g 
(\int_{\gamma_\nu}dz_i \int_{\gamma_{g+\nu}} {\mathfrak X}+\mu_q
-\int_{\gamma_{g+\nu}}dz_i \int_{\gamma_{\nu}} {\mathfrak X}+\mu_q) 
\right)_{i=1,\ldots,g}\in \Jac \bar{X}.}
\]

\subsection{A Higher Reciprocity Law}

Generally for functions $F,G:\pi_1(X,p) \rightarrow
\cpxm$ we introduce: $\Period(F;\,G):=
\sum_{\nu=1}^g \big(F({{\gamma}_\nu})G({{\gamma}_{g+\nu}}) 
-F({{\gamma}_{g+\nu}}) G({{\gamma}_{\nu}})\big)$. For instance,
Riemann's classical 
period relation reads: $\Period(\int dz_i;\int dz_j)=0$. 
With this notation we can state a higher reciprocity law.  

\begin{theorem} \label{Elefant}
For any holomorphic 1-form $\omega$ on $\bar{X}$ holds:
\begin{multline*}
\sum_{\nu=1}^g \{ \int_{\gamma_\nu}\omega \int_{\gamma_{g+\nu}}
{\mathfrak X} + \mu_q
- \int_{\gamma_{g+\nu}}\omega \int_{\gamma_{\nu}}{\mathfrak X} + \mu_q \} 
= 2g \int_p^q \omega \\ 
+\sum_{j,k=1}^g  a_{jk}\left\{
\Period(\int \omega;\int dz_j\int d\bar{z}_k) 
+2\Period(\int \omega \int d\bar{z}_k;\int dz_j)\right. 
\left.  -2\Period(\int \omega dz_j;\int d\bar{z}_k)     
\right\}.
\end{multline*}
\end{theorem}               
                                     
\subsubsection{Observation} \label{Zwerg Nase}

The proof of the {\it higher 
reciprocity law} in Theorem \ref{Elefant} and also later the proof
of the {\it higher period relation} of Theorem \ref{Tante} are direct
generalizations of Riemann's bilinear relations as they are proved 
in \cite{Chen} or \cite{Gunning-quadratic-periods}. Similar considerations
like the following can also be found in \cite{Poor-Yuen}. 

Let $c_i:=(\gamma_i-1)$ and $d:=(\delta-1)$ denote the elements in $J$ 
corresponding to $\gamma_i$ and $\delta$ in $\pi_1(X,p)$. 
If we interpret relation \eqref{Hannelore} in
$\itgm\pi_1(X,p)$ modulo $J^4$, we obtain 
\begin{multline} \label{Gunter}
\sum\nolimits_{\nu=1}^g \{ c_\nu c_{g+\nu} - c_{g+\nu} c_\nu \\
+ (c_{g+\nu} c_{\nu} c_{g+\nu} - c_{\nu} c_{g+\nu} c_\nu)
- (c_\nu c_{g+\nu} c_{g+\nu} - c_{g+\nu} c_\nu c_{\nu}) \}
 \equiv d  \mod{J^4}.
\end{multline}                  
When the linear extension of a homotopy functional 
$F:\pi_1(X,p)\rightarrow \cpxm$
to $\itgm\pi_1(X,p)$ satisfies $F(J^4)=0$, then it has to respect this relation.
For instance iterated integrals of length $\le 3$, which are homotopy 
functionals, are examples for such $F$.

\noindent
{\bf Proof of Theorem \ref{Elefant}:}
Use that for a closed path $\alpha$ holds:
$\int_{\alpha}d\bar{z}_j dz_k + \int_{\alpha} dz_k d\bar{z}_j = 
\int_{\alpha}d\bar{z}_j  \int_{\alpha} dz_k$ to prove that the left 
hand side of the equation in Theorem \ref{Elefant} equals
\[
\Period(\int \omega;\int\sum_{j,k=1}^g 2a_{jk} dz_j d\bar{z}_k + \mu_q)
-\Period(\int \omega;\sum_{j,k=1}^g a_{jk}\int dz_j\int d\bar{z}_k ).
\]   

Note that $\int I:=\int\sum_{j,k=1}^g 2a_{jk} \omega dz_j d\bar{z}_k 
+ \omega\mu_q$ is
a homotopy functional, so its values on both sides of \eqref{Gunter}
coincide. Recall that for 1-forms $\varphi,\,\psi,\,\chi$ and closed paths 
$\alpha$, $\beta$ with $a=(\alpha-1)$, $b=(\beta-1)$ and 
$ab=(\alpha\beta-\alpha-\beta+1)$ holds: 
$\int_{ab} \varphi\psi\chi=\int_{a} \varphi\int_{b}\psi\chi
+\int_{a} \varphi\psi\int_b \chi$. Using this rule, a direct computation 
shows that the value of $\int I$ on the left hand side of 
relation \eqref{Gunter} takes the value:
\begin{multline*}
\Pi(\int\omega;\int I) \\
+\sum_{j,k=1}^g 2a_{jk}\left\{
 \Period(\int \omega dz_j;\int d\bar{z}_k)
-\Period(\int \omega \int d\bar{z}_k;\int dz_j)
-\Period(\int \omega; \int dz_j\int d\bar{z}_k)
\right\}.
\end{multline*}                 

According to our observation \ref{Zwerg Nase} this has to be equal
to the value of the homotopy functional $\int I$ applied to the right 
hand side of \eqref{Gunter}. We compute this value as follows.  
From $\wedge{\mathfrak X}+d\mu_q=0$ we can determine the shape
of $\mu_q$. Using Stokes' theorem,
a standard argument shows that there is a simply connected holomorphic
coordinate plot
$(U,z)$ on $\bar{X}$ containing $q$ and all of a representing path for
$\delta\in\pi_1(X,p)$ such that on $U$ we may write
$\mu_q=\frac{2g}{2\pi i}\; \frac{dz}{z} +\varphi$, 
where $\varphi$ is a smooth (non-closed) 1-form in $E^1(U)$.    
Since this representative of $\delta$ is 0-homotopic in $U$, the 
homotopy functional $\sum_{j,k=1}^g 2a_{jk} \int \omega dz_j d\bar{z}_k      
+\omega \varphi$ vanishes on it. It remains: 
$\int_\delta \omega (\frac{2g}{2\pi i}\frac{dz}{z})=2g\int_p^q\omega$.
Putting all ingredients together provides the proof. \qed

\subsection{A Higher Period Relation}

Recall that we chose $\Omega=(I,Z)$ with symmetric $Z$ having positive 
imaginary part and with $A=(\bar{Z}-Z)^{-1}$. 
Define for $i=1,\ldots,g$ the $g\times g$-matrices
\[
I^i_1:=\left(\;\int_{c_\nu} dz_i dz_j \right)_{\nu,j}\text{ and }
I^i_2:=\left(\;\int_{c_{g+\nu}} dz_i dz_j \right)_{\nu,j}
 \in \Mat(g\times g;\cpxm).
\]
Then we define the following two vectors with entries in 
$\Mat(g\times g;\cpxm)$
\[
{I_1=
\begin{pmatrix}I_1^1 \\ \vdots \\ I_1^g\end{pmatrix},\quad I_2=
\begin{pmatrix}I_2^1 \\ \vdots \\ I_2^g\end{pmatrix}}\quad
{\in \Mat\big(g\times 1; \Mat(g\times g)\big)}.
\]
For some matrix $M$, denote by $\tr M$ the {\it trace} of $M$ and by $\diag M$
its {\it diagonal}. For a vector of matrices let {\it the trace} of this
vector be the vector of the traces of the matrices. 
The following theorem is the announced {\it higher period relation}. 

\begin{theorem} \label{Tante}
\begin{multline*}
\left( 2 \tr (I_2 A) - 2 \tr (I_1 AZ )\right)
+ \left(\diag (ZAZ) - Z \diag (AZ)\right)  \\
+ \left(\diag (ZA) - Z \diag (A)\right)
+ \left(\diag (AZ) - Z \diag (ZA)\right)  \equiv 0
\mod{(I,Z)\itgm^{2g}}. 
\end{multline*}
\end{theorem}

\noindent
{\bf Proof:}
Apply the homotopy functional $\sum_{j,k=1}^g a_{jk} 
\int dz_j\,dz_i\,dz_k$ to \eqref{Gunter}. \qed

\vspace{1ex}
With the above notation, we use this higher period relation to continue 
our computation of the extension $k_{pq}$. After
Theorem \ref{Elefant} it makes sense to speak of $k_{pq}$; we have
$\Psi(k_{pq})=2g(q-p)+\Psi(k_{pp})$. 
For $u\circ\Psi(k_{pp})\in \cpxm^g\big/\Omega\itgm^{2g}$ we had the 
following expression: \qquad  \qquad $u\circ\Psi(k_{pp})=$
\[
\diag (Z A\bar{Z}) - Z \diag (A)
+ 2 \diag (ZA) - 2 Z \diag (A\bar{Z})
- 2 \tr (I_1 A \bar{Z}) + 2 \tr (I_2 A)
\]
Transform this expression such that it only contains (iterated)
integrals over {\it holomorphic} forms.
Observe
$\diag (ZA\bar{Z})=\diag (Z(\bar{Z}-Z)^{-1}(\bar{Z}-Z))
+\diag (ZA{Z})
=\diag (Z)+\diag (ZA{Z})$
and similarly
$2 Z \diag (A\bar{Z})\equiv 2 Z\diag (AZ)
\mod{(I,Z)\itgm^{2g}}$
and $2 \tr (I_1 A\bar{Z})=2 \tr (I_1)+2 \tr (I_1 A{Z})$.
Using these identities we continue
\begin{eqnarray*}
u\circ\Psi(k_{pp})
& \equiv & \diag (Z) + \diag (ZAZ) - Z\diag (A)\\
& + & 2 \diag (ZA) - 2 Z \diag (A{Z}) \\
& - & 2 \tr (I_1) - 2 \tr (I_1 A {Z})
+ 2 \tr (I_2 A) \mod{(I,Z)\itgm^{2g}}.
\end{eqnarray*}
Notice: $\diag (ZA) - \diag(AZ)=0$. When we apply Theorem \ref{Tante}, we 
finally get: $u\circ\Psi(k_{pp}) \equiv \diag Z - 2 \tr (I_1) 
\mod{(I,Z)\itgm^{2g}}$. Writing this out, we find 
$u\circ\Psi(k_{pp})\equiv 2\kappa_p\mod{(I,Z)\itgm^{2g}}$ 
by virtue of formula 
\eqref{Rumpelstielzchen}. This is the proof of
Theorem \ref{Hans guck in die Luft}.   

\section{A Pointed Torelli Theorem for  Punctured Riemann 
Surfaces} \label{Gretel}

Here we want to show that the extension $w_{pq}$ 
or respectively $m_{pq}$ determines $p$ and $q$. For this we do not need
Section \ref{Haensel}. Finally we will combine this with
results of Hain and Pulte \cite{Hain-the-geometry}, \cite{Pulte}, 
which we briefly sketch first.

\subsection{The Pointed Torelli Theorem}

The pointed Torelli Theorem of Hain and Pulte is based on the following.
\begin{theorem}[Hain, Pulte] \label{injection}
The map from 
$\Pic^0 \bar{X}$ to $\Extmhs (K;H^1)$
which maps  $(p-p')$ to $m_p-m_{p'}$ is well-defined and injective.
\end{theorem}

We write $(\bar{X},p)\cong (\bar{X},p')$ if there is an automorphism 
$\phi:\bar{X}\rightarrow \bar{X}$
that maps $p$ to $p'$. For a point $p$ on $\bar{X}$ we define {\it the
set of alternatives for $p$} as 
\[
a_{\bar{X}}(p):=\{p\}\cup\left\{p'\in \bar{X}| 
m_{p'}=-m_p \text{ and } (\bar{X},p)\not\cong (\bar{X},p')\right\}
\]
The following is a consequence of Theorem \ref{injection}.
Let us give a short proof of it.

\begin{corollary}
$a_{\bar{X}}(p)$ consists of at most two points. Up to automorphism of 
$\bar{X}$,
there cannot be more than one pair of different points $\{p,p'\}$ 
on $\bar{X}$ such that $a_{\bar{X}}(p)=\{p,p'\}=a_{\bar{X}}(p')$. 
\end{corollary}

\noindent
{\bf Proof:} 
The first assertion is an obvious consequence of \ref{injection}. To 
prove the second assertion, assume that $\tilde{p}$ and $\tilde{p}'$ is 
another such pair with $a_{\bar{X}}(\tilde{p})=\{\tilde{p},\tilde{p}'\}
=a_{\bar{X}}(\tilde{p})$. Then by \ref{injection}, the divisors $p+p'=
\tilde{p}+\tilde{p}'$ are linearly equivalent. It follows that 
either $\{p,p'\}=
\{\tilde{p},\tilde{p}'\}$ or $\bar{X}$ is hyperelliptic and the hyperelliptic
involution maps $p$ to $p'$ and $q$ to $q'$, which contradicts
the assumptions on $p$, $p'$ and $q$, $q'$. \qed

Together with the classical Torelli theorem, Hain and Pulte used 
Theorem \ref{injection} to prove the following {\it pointed Torelli theorem}.
For a pointed compact Riemann surface $(\bar{Z},z_0)$ denote by 
$J_{z_0}(\bar{Z})$ the augmentation ideal in $\itgm\pi_1(Z,z_0)$.

\begin{theorem}[Hain, Pulte]  \label{pointed Torelli theorem}
Suppose that $(\bar{X},p)$ and $(\bar{Y},r)$ are two pointed compact Riemann
surfaces. If there is a ring homomorphism
\[
\itgm \pi_1(\bar{X},p) \Big/ J_p(\bar{X})^3
\xrightarrow{\cong}
\itgm \pi_1(\bar{Y},r) \Big/ J_r(\bar{Y})^3
\]
which induces an isomorphism of MHSs, then there is an isomorphism
$f:\bar{X}\rightarrow \bar{Y}$ with $f(p)\in a_{\bar{Y}}(r)$.
\end{theorem}

\begin{remark} \rm
As far as the author knows, still no example is known of a pointed 
compact Riemann surface $(\bar{X},p)$ with $|a_{\bar{X}}(p)|=2$. 
M.~Pulte \cite{Pulte} has shown that such an $(\bar{X},p)$ with 
$|a_{\bar{X}}(p)|=2$ 
must have zero harmonic volume. B.~Harris \cite{Bruno-Harris} 
proved that a generic smooth projective complex curve has non zero 
harmonic volume. Moreover, Pulte showed (loc.~cit.) that, if there
are two points $p$, $p'$ with $a_{\bar{X}}(p)=\{p,p'\}=a_{\bar{X}}(p')$, then 
$(g-1)(p+p')-K=0\in \Pic^0 \bar{X}$, where $K$ is the canonical divisor. 
For pointed hyperelliptic curves $(\bar{X},p)$ always holds: 
$a_{\bar{X}}(p)=\{p\}$,
since here $m_p=-m_{p'}$ implies $(\bar{X},p)\cong(\bar{X},p')$ 
by the hyperelliptic involution. 
\end{remark}

\subsection{A Punctured Pointed Torelli Theorem}

The following Theorem will follow directly from 
Lemma \ref{Tischlein deck dich!}, which we prove at the end of this section.  
\begin{theorem}             \label{double injection}
For all $p\in \bar{X}$, the map 
$\Pic^0 \bar{X} \rightarrow \Extmhs (H^1\otimes H^1;H^1)$
which maps  $(q-q')$ to $m_{pq}-m_{pq'}$ is well-defined and injective.
\end{theorem}

Combining Theorem \ref{double injection} with 
the results of Hain and Pulte we find.
\begin{proposition}   \label{Josef Beuys}
The map from 
$(\bar{X}\times \bar{X})\setminus \Delta$ to $\Extmhs (H^1\otimes H^1;H^1)$ given by
$(p,q)\mapsto m_{pq}$
is well-defined, extends to the diagonal $\Delta$ and is injective.
\end{proposition}

\noindent
{\bf Proof of \ref{Josef Beuys}:}
Note that the map of complex tori $\Extmhs(H^1\otimes H^1;H^1)\rightarrow
\Extmhs(K\oplus Q;H^1)$ is a covering map, since
$\Hom(H^1\otimes H^1;H^1)_\itgm\hookrightarrow \Hom(K\oplus Q;H^1)_\itgm$.
Moreover, we have the commutative diagram
\[
\begin{CD}
(X\times X)\setminus \Delta @>\tilde{\varphi}>> \Extmhs(H^1\otimes H^1;H^1) \\
@VVV @VV\begin{matrix}\text{\scriptsize covering}\\
\text{\scriptsize map}\end{matrix}V \\
X\times X @>\varphi>> \Extmhs(K;H^1)\oplus \Pic^0 \bar{X} \\
(p,q) & \longmapsto & \big( m_p,\, (2g\, q -2p-K)\big).
\end{CD}
\]

The map $\varphi$ is continuous ($m_p$ is -- in a coordinate system --
an expression of iterated integrals over paths with basepoint $p$). As
the map $\tilde{\varphi}(p,q)=m_{pq}$ is a lifting of $\varphi$,
we see that $\tilde{\varphi}$ is continuous too. The fact that the
map $m_{pq}\mapsto (m_p,\,k_{pq})$ is a covering map tells us moreover
that we may extend $\tilde{\varphi}$ to the diagonal $\Delta$.
Now by \ref{injection}, the result of Hain and Pulte, the extension
$m_{pq}$ determines $p$. By Theorem \ref{double injection} it determines
also $q$. \qed

Pulling back the intersection form 
$H_1(\bar{X},\itgm)\otimes H_1(\bar{X},\itgm)\rightarrow \itgm$ along
the natural isomorphism $J\big/ J^2\xrightarrow{\cong}\bar{J}\big/ \bar{J}^2$
induces a polarization on $\Gr_{-1}^W(J\big/ J^3)=H^1(X)$. We can also
put a polarization on $\Gr_{-2}^W(J\big/ J^3)=J^2\big/ J^3\cong
J\big/ J^2\otimes J\big/ J^2=
\Gr_{-1}^W(J\big/ J^3)\otimes\Gr_{-1}^W(J\big/ J^3)$, by taking the 
tensor product of the polarized Hodge structure $H^1(X)$ in the 
category of polarized Hodge structures. In that sense, $J\big/ J^3$ becomes 
a {\it graded polarized MHS}, i.~e.~each $\Gr_{l}^W$ is a polarized Hodge
structure.

For points $p$ and $q$ on $\bar{X}$ we define
\[
A_{\bar{X}}(p,q):=\{(p,q)\}\cup\left\{(p',q')\in \bar{X}\times \bar{X} 
\; \left|  \begin{array}{l}
m_{p'q'}=-m_{pq} \quad\text{ and } \\ 
\Big(\bar{X}\setminus \{q\},p\Big)\not\cong 
\Big(\bar{X}\setminus \{q'\},p'\Big)\end{array} \right.\right\}.
\]
The following is then a consequence of Proposition \ref{Josef Beuys}. 

\begin{corollary}
$A_{\bar{X}}(p,q)$ consists of at most two elements. \qed
\end{corollary}

Our results lead to the following {\it punctured pointed Torelli theorem}.
\begin{theorem}  \label{two-pointed Torelli theorem}
Suppose that $(\bar{X}\setminus\{q\},p)$ and $(\bar{Y}\setminus\{s\},r)$ 
are two punctured
compact Riemann surfaces with basepoint. If there is a ring isomorphism 
\[
\itgm \pi_1(\bar{X}\setminus\{q\},p) \Big/ J_p(\bar{X}\setminus\{q\})^3
\xrightarrow{\cong}
\itgm \pi_1(\bar{Y}\setminus\{s\},r) \Big/ J_r(\bar{Y}\setminus\{s\})^3,
\]
which induces an isomorphism of graded polarized MHSs, then there is
a biholomorphism $f:\bar{X}\rightarrow \bar{Y}$ with 
$(f(p),f(q))\in A_{\bar{Y}}(r,s)$.
\end{theorem}

\noindent
{\bf Proof of \ref{two-pointed Torelli theorem}:}
The proof goes along the lines of the proof of the pointed Torelli theorem in
\cite{Pulte} and \cite{Hain-the-geometry}. 
Let $J_{pq}=J_p(\bar{X}\setminus\{q\})$ and $J_{rs}
=J_r(\bar{Y}\setminus\{s\})$.
We have an isomorphism of MHSs, $\lambda:J_{pq}\big/J_{pq}^3
\xrightarrow{\cong} J_{rs}\big/J_{rs}^3$ and in particular, $\lambda$
induces an isomorphism of polarized Hodge structures
\[
\lambda^*:
H^1(\bar{Y})=W_1\left(J_{rs}\big/J_{rs}^3\right)^* \rightarrow
W_1\left(J_{pq}\big/J_{pq}^3\right)^*= H^1(\bar{X}).
\]
By the classical Torelli theorem
(cf. for instance \cite{Martens}) we know that there is a biholomorphism
$f:\bar{X}\rightarrow \bar{Y}$ such that $f^*:H^1(\bar{Y})
\rightarrow H^1(\bar{X})$ is $\pm \lambda^*$.
Since $\lambda$ respects the {\it ring} structure, the by $\lambda$ induced map
$\left(J_{rs}^2\big/ J_{rs}^3\right)^*\rightarrow
\left(J_{pq}^2\big/ J_{pq}^3\right)^*$
is determined by $\lambda^*:H^1(\bar{Y})\rightarrow H^1(\bar{X})$ and hence,
\[
f^*: \left(J_{rs}^2\big/ J_{rs}^3\right)^*=H^1(\bar{Y})\otimes H^1(\bar{Y})
\rightarrow
H^1(\bar{X})\otimes H^1(\bar{X})=\left(J_{pq}^2\big/ J_{pq}^3\right)^*
\]
is equal to $\lambda^*\otimes\lambda^*$.
Without loss of generality, we may therefore assume that 
$(\bar{Y}\setminus\{s\},r) = (\bar{X}\setminus\{q'\},p')$ for 
two points $p'$ and $q'$ in $\bar{X}$ and that
the following diagram commutes
\[
\begin{CD}
0 @>>> H^1 @>>> \left( J_{pq}\big/ J_{pq}^3\right)^* @>>>
H^1\otimes H^1 @>>> 0 \\
& & @V{\pm id}VV @VV{\lambda^*}V @VV{id}V \\
0 @>>> H^1 @>>> \left( J_{p'q'}\big/ J_{p'q'}^3\right)^* @>>>
H^1\otimes H^1 @>>> 0 .
\end{CD}
\]
It follows that $m_{pq}=\pm m_{p'q'}$. This means that there either
is an automorphism $\phi:(\bar{X}\setminus\{q\},p)\rightarrow 
(\bar{X}\setminus\{q'\},p')$
or $A_{\bar{X}}(p,q)=\{(p,q);\;(p',q')\}=
A_{\bar{X}}(p',q')$. In both cases, the identity
map is the map with the desired properties. \qed

\subsection{A Technical Lemma}

Theorem \ref{double injection} is a consequence of the following
\begin{lemma} \label{Tischlein deck dich!}
For any element $\sum_i (q_i-q_i')\in \Pic^0 \bar{X}$ holds:
\[
\sum\nolimits_i (q_i-q_i')=0 \in \Pic^0 \bar{X} \; \Leftrightarrow\; 
\sum\nolimits_i (m_{pq_i}-m_{pq_i'})=0 \in \Extmhs \big((H^1)^{\otimes 2};H^1\big).
\]
\end{lemma}

\noindent
{\bf Proof:}
Consider the isomorphism (cf.~\cite{Carlson})
from $\Extmhs \big(H^1\otimes H^1;H^1\big)$ to 
\[
{\Hom (H^1\otimes H^1; H^1)_\cpxm}\Big/ 
(F^0 \Hom (H^1\otimes H^1; H^1)_\cpxm+\Hom (H^1\otimes H^1; H^1)_\itgm),
\]
where the image of an extension $m_{pq}$ is $[\phi_{pq}]$ for a certain
$\phi_{pq}\in\Hom (H^1\otimes H^1; H^1)_\cpxm$. 
On an element $[\varphi]\otimes [\psi]\in H^1\otimes H^1$, the homomorphism
$\phi_{pq}$ has the following property. 
There is a 
$\eta_q\in F^1 E^1(X\log q)$ such that $\varphi\wedge\psi + d\eta_q=0$ and
$\phi_{pq}\left([\varphi]\otimes [\psi]\right)= \sum_{j=1}^{2g}
(\int_{\gamma_j} \varphi\psi+\mu_q) [dx_j]$.
If $[\varphi]\otimes [\psi]\in K$ then $\eta_q$ can be chosen in 
$F^1 E^1(X)$ and does not depend on $q$, which shows that
$\left(\phi_{pq}-\phi_{pq'}\right)$ is zero on $K$.
Therefore
it is determined by its value on one element of $\left(H^1\otimes H^1\right)
\setminus K$; for instance on $[dx_1]\otimes [dx_{g+1}]$. 

Given a divisor $D=\sum_i (q_i-q_i')$ define the homomorphism 
$\Phi_D:= \sum_i \left(\phi_{pq_i}-\phi_{pq_i'}\right):H^1\otimes H^1 
\rightarrow H^1$.
We will derive a series of equivalences. First, we have: 
\[
\begin{split}
 &\sum\nolimits_i \left(m_{pq_i}-m_{pq_i'}\right)= 0 \in 
\Extmhs \big(H^1\otimes H^1;H^1\big) \\
\Leftrightarrow \quad & \Phi_D \in F^0 \Hom \big(H^1\otimes H^1;H^1\big)_\cpxm +
\Hom \big(H^1\otimes H^1;H^1\big)_\itgm .
\end{split}
\]

Now let ${\mathbf w}\in H^{0,1}\otimes H^{0,1}$ be such that 
$[dx_1]\otimes [dx_{g+1}]-{\mathbf w} \in F^1 (H^1\otimes H^1)= 
H^{1,0}\otimes H^{1}+ H^{1}\otimes H^{1,0}$. Note that $H^{0,1}\otimes H^{0,1}
\subset K$ and hence $\Phi_D({\mathbf w})=0$. Moreover $H^{1,0}\otimes H^{1,0}
\subset K$ and $\Phi_D\left(H^{1,0}\otimes H^{1,0}\right)=0$.
Therefore, we may continue the series of equivalences by:
\[
\Leftrightarrow \Phi_D\left( [dx_1]\otimes [dx_{g+1}]-{\mathbf w} \right)\in 
H^{1,0}+H^1_\itgm 
\Leftrightarrow \Phi_D\left( [dx_1]\otimes [dx_{g+1}] \right)\in 
H^{1,0}+H^1_\itgm. 
\]
Let $\eta_{q_i}\in F^1 E^1(X\log q_i)$ and $\eta_{q_i'}\in
F^1 E^1(X\log q_i')$ be such that 
$dx_1\wedge dx_{g+1} + d\eta_{q_i}=0$ 
and $dx_1\wedge dx_{g+1} + d\eta_{q_i'}=0$.
Note that this implies
$\Res_{q_i} \eta_{q_i}=\frac{1}{2\pi i}=\Res_{q_i'} \eta_{q_i'}$.
Then a direct computation shows that we may go on:
\[
\begin{split}
\Leftrightarrow \quad & \sum\nolimits_{j=1}^{2g} \sum\nolimits_i
(\int_{\gamma_j}\mu_{q_i}-\mu_{q_i'})[dx_j]\in
H^{1,0}+H^1_\itgm \\
\Leftrightarrow \quad & \left( \Period \big(\int dz_\nu;\; 
\int (\mu_{q_i}-\mu_{q_i'})\big)\right)_\nu \equiv 0 \mod{\Omega \itgm^{2g}}.
\end{split}
\]
By the reciprocity law for differentials of the third kind 
(cf.~\cite{Griffiths-Harris}), we find as $(\mu_{q_i}-\mu_{q_i'})$ is 
meromorphic with simple poles:
$\Leftrightarrow  \sum_i(q_i-q_i') =0 \in \Pic^0 \bar{X}$. 
That proves the lemma. \qed

% LITERATURLISTE


\begin{thebibliography}{Mum83}

\bibitem[Car80]{Carlson}
J.~A. Carlson.
\newblock Extensions of mixed {Hodge} structures.
\newblock In {\em {Journ{\'e}es de G{\'e}om{\'e}trie Alg{\'e}brique d'Angers}},
  pages 107--128, Alphen aan den Rijn, 1980. Sijthoff and Noordhoff.

\bibitem[Che76]{Chen-RBC}
K.~T. Chen.
\newblock Reduced bar constructions on de {Rham} complexes.
\newblock In A.~Heller and (eds.) M.~Tierny, editors, {\em Algebra, Topology
  and Category Theory}, pages 19--32, New York, 1976. Academic Press.

\bibitem[Che77]{Chen}
K.~T. Chen.
\newblock Iterated path integrals.
\newblock {\em Bull.~Amer.~Math.~Soc.}, 83:831--879, 1977.

\bibitem[Fay73]{Fay}
J.~D. Fay.
\newblock {\em Theta Functions on {Riemann} Surfaces}.
\newblock Number 352 in Lecture Notes in Math. Springer Verlag, Berlin, 1973.

\bibitem[GH78]{Griffiths-Harris}
P.~Griffiths and J.~Harris.
\newblock {\em Principles of Algebraic Geometry}.
\newblock John Wiley \& Sons, Inc., USA, 1978.

\bibitem[Gun69]{Gunning-quadratic-periods}
R.~C. Gunning.
\newblock Quadratic periods of hyperelliptic abelian integrals.
\newblock In {\em Problems in Analysis}, pages 239--247, New Jersey, 1969.
  Princeton University Press.

\bibitem[Hai87a]{Hain-de-Rham-homotopy}
R.~M. Hain.
\newblock The de {Rham} homotopy theory of complex algebraic varieties {I}.
\newblock {\em K-Theory}, 1:271--324, 1987.

\bibitem[Hai87b]{Hain-the-geometry}
R.~M. Hain.
\newblock The geometry of the mixed {Hodge} structure on the fundamental group.
\newblock {\em Algebraic Geometry, Bowdoin, PSPM}, 46:247--282, 1987.

\bibitem[Har83]{Bruno-Harris}
B.~Harris.
\newblock Harmonic volumes.
\newblock {\em Acta Math.}, 150(1--2):91--123, 1983.

\bibitem[HL97]{Hain-Looijenga}
R.~M. Hain and E.~J.~N. Looijenga.
\newblock Mapping class groups and moduli spaces of curves.
\newblock {\em Proc. Symp. Pure Math.}, 62(pt.II):97--142, 1997.
\newblock Algebraic geometry, Santa Cruz.

\bibitem[Jab86]{Jablow}
E.~R. Jablow.
\newblock Quadratic vector classes on {Riemann} surfaces.
\newblock {\em Duke Math.~J.}, 53(1):221--232, 1986.

\bibitem[Lan02]{Landfriedt}
E.~Landfriedt.
\newblock {\em {Thetafunktionen und hyperelliptische Funktionen}}.
\newblock G{\"o}schensche Verlagshandlung, Leipzig, 1902.

\bibitem[Lew64]{Lewittes}
J.~Lewittes.
\newblock Riemann surfaces and the theta function.
\newblock {\em Acta Math.~}, 111:37--61, 1964.

\bibitem[Mar63]{Martens}
H.~H. Martens.
\newblock A new proof of {Torelli's} theorem.
\newblock {\em Annals of Math.}, 78(1):107--111, 1963.

\bibitem[Mum83]{Tata-on-theta-I}
D.~Mumford.
\newblock {\em Tata Lectures on Theta I}.
\newblock Number~28 in Progress in Math. Birkh{\"a}user, Boston, 1983.

\bibitem[Pul88]{Pulte}
M.~Pulte.
\newblock The fundamental group of a {Riemann} surface: mixed {Hodge}
  structures and algebraic cycles.
\newblock {\em Duke Math.~J.}, 57(3):721--760, 1988.

\bibitem[PY96]{Poor-Yuen}
C.~Poor and D.~S. Yuen.
\newblock Relations on the period mapping giving extensions of mixed {Hodge}
  structures on compact {Riemann} surfaces.
\newblock {\em Geometria Dedicata}, 59:243--291, 1996.

\bibitem[Rie92]{Riemann}
B.~Riemann.
\newblock {\em {Bernhard Riemann's gesammelte mathematische Werke}}.
\newblock B.~G.~Teubner, Leipzig, 1892.

\end{thebibliography}
\end{document}